\documentclass[final,times,twocolumn]{elsarticle}
\usepackage[table, svgnames, dvipsnames]{xcolor}
\usepackage{lineno,hyperref}
\modulolinenumbers[5]
\usepackage{natbib}
\usepackage{tabularx}
\usepackage{comment}
\usepackage{color}
\usepackage{multicol,tikz}
\usepackage{xspace,amssymb,amsmath,epsfig}
\usepackage{multirow}
\usepackage[T1]{fontenc}
\usepackage{spverbatim}
\usepackage{fancyvrb}
\usepackage{fvextra}
\usepackage{graphicx,multirow}
\usepackage{framed,booktabs}

\newcommand{\LDLT}{$\textrm{LDL}^{\textrm{T}}$}

\journal{Parallel Computing}

\begin{document}
\begin{frontmatter}
\title{Linear solvers for power grid optimization problems:
    \\ a review of GPU-accelerated linear solvers}

\author[nrel,pnnl]{Kasia \'{S}wirydowicz}
\author[su]{Eric Darve}
\author[nrel]{Wesley Jones}
\author[nrel]{Jonathan Maack}
\author[su]{Shaked Regev}
\author[su]{Michael A. Saunders}
\author[nrel]{Stephen J. Thomas}
\author[pnnl]{Slaven Pele\v{s}}

\address[pnnl]{Pacific Northwest National Laboratory}
\address[nrel]{National Renewable Energy Laboratory}
\address[su]{Stanford University}

\date{\today}

\begin{abstract}
The linear equations that arise in interior methods for constrained optimization are sparse symmetric indefinite, and they become extremely ill-conditioned as the interior method converges. These linear systems present a challenge for existing solver frameworks based on sparse LU or \LDLT{} decompositions. We benchmark five well known direct linear solver packages on CPU- and GPU-based hardware, using matrices extracted from power grid optimization problems. The achieved solution accuracy varies greatly among the packages. None of the tested packages delivers significant GPU acceleration for our test cases. For completeness of the comparison we include results for MA57, which is one of the most efficient and reliable CPU solvers for this class of problem.
\end{abstract}

\end{frontmatter}

\section{Introduction}

Interior-point methods \cite{waecther_05_IPOPT,waecther_05_IPOPT2} provide an important tool for nonlinear optimization and are used widely for engineering design \cite{biegler2013survey} and discovery from experimental data \cite{duan2010three}. Most of the computational cost for interior methods comes from solving underlying linear equations. Indeed, interior methods came to prominence with the emergence of robust sparse linear solving techniques. For a general nonconvex optimization problem, the linear systems are sparse symmetric indefinite and typically ill-conditioned. Furthermore, many implementations require the linear solver to provide matrix inertia information (the number of positive, zero, and negative eigenvalues). Relatively few presently available linear solvers meet the requirements of interior methods. 

We review the current state-of-the-art, with particular emphasis on solvers that can run on hardware accelerators such as GPUs. Our investigation is motivated by power grid optimization work within the  US Department of Energy’s Exascale Computing Project. In a large-scale security-constrained optimal power flow analysis, several coupled optimal power flow problems are solved simultaneously \cite{chakrabarti2014security,petra2014real}. Each coupled unit runs its own interior method computation. Because more than 95\% of the computational power on new and emerging computational platforms is typically in accelerators, it is desirable to run each interior method optimization, including linear solves, on those devices. We identify technical gaps in this area and suggest possible paths to bridge those gaps.

The paper is organized as follows: In Section \ref{sec:problem_statement}, we provide a brief description of the algorithms employed in interior point method and we formulate the problem; we also present an overview of the available literature, and place the problem in a broader context. In section \ref{sec:linear_solvers}, we discuss the main features of the linear solver packages. In section \ref{sec:tests}, we describe the test cases. In section \ref{sec:results}, we present test results for each software package. These results are summarized in Section \ref{sec:summary}. In section \ref{sec:conclusion}, we state conclusions and provide ideas for future work.

\section{Problem statement} \label{sec:problem_statement}

We are concerned with solving large linear systems
\begin{equation}
Ax =b,
\label{eq:mathproblemstatement}
\end{equation}
where $A$ is a symmetric indefinite $n \times n$ sparse matrix, $b$ is an $n \times 1$ right-hand-side (RHS) vector, and $x$ is the $n \times 1$ solution vector. The  \emph{residual vector} is $r =b-Ax$.

We generate our test problems for linear solvers using the Ipopt optimization package \cite{wachter2006implementation}, which is considered to be the ``gold standard'' for freely available interior method software. We instrumented Ipopt to save $A$ and $b$, which are handed off to the linear solver every time $A$ is updated. Each series of test cases corresponds to a successful optimization solver run, i.e., Ipopt converged to the optimal solution. Ipopt used MA57 \cite{Duff2004} as the linear solver in these computations.

The interior method finds the minimum of a scalar objective function given equality and inequality constraints. The method extends the minimization problem with barrier functions in order to enforce inequality constraints and bounds on the variables. The extended problem is solved using a variant of Newton's method. A continuation is applied to reduce barrier parameters gradually and bring the extended problem as close to the original problem as possible. A nonlinear problem is (re)solved at each continuation step. The problem becomes singular when barrier parameters are zero, so the continuation needs to stop when the extended problem solution is close enough to the solution of the original problem, but before the problem becomes too ill-conditioned. This is why having a robust linear solver is critical for successful implementation of the interior method. For more details, the reader is referred to \cite{wachter2006implementation}.

\begin{figure}
    \centering
    \includegraphics[width=\columnwidth]{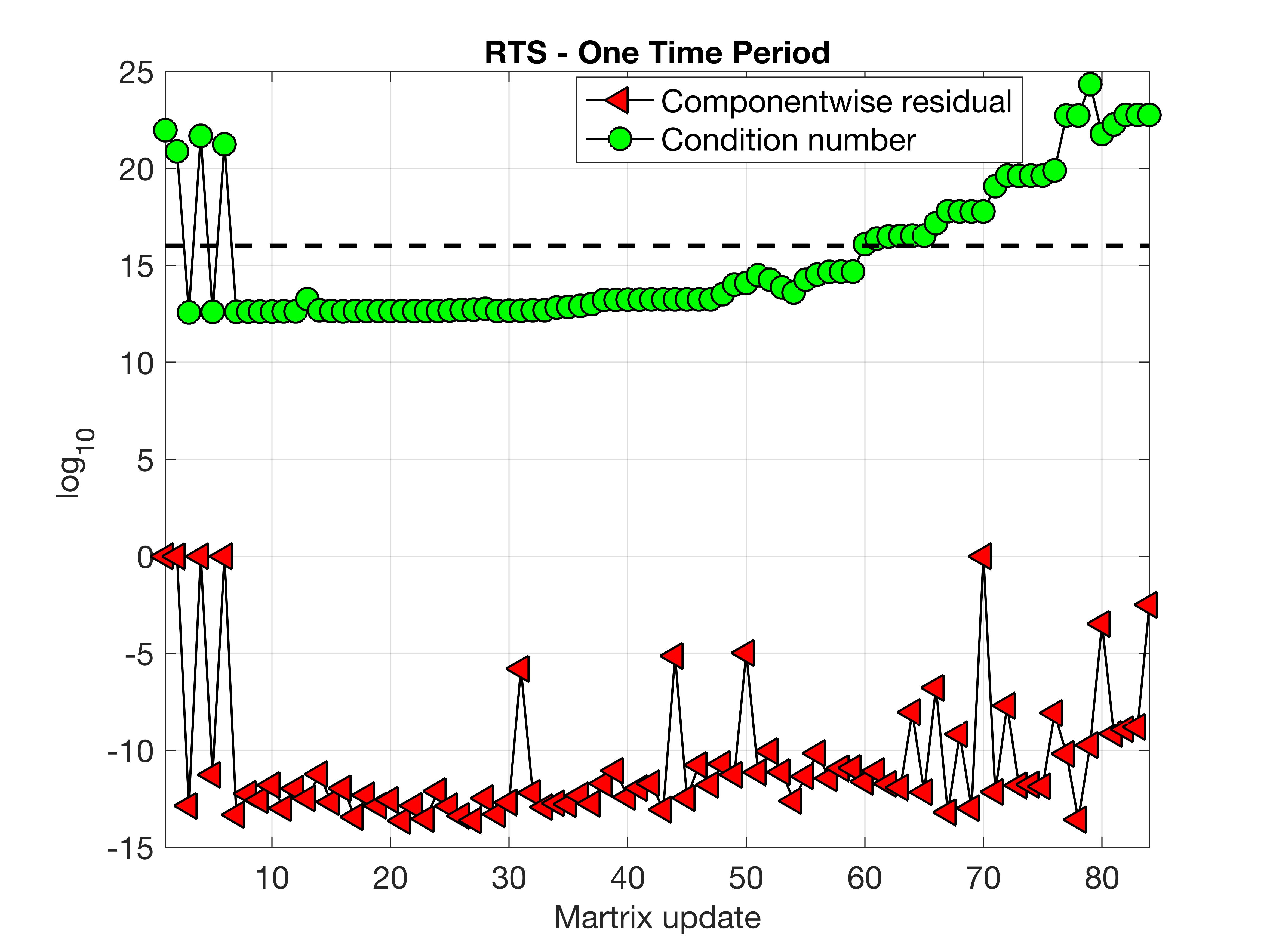}
    \caption{Condition number and relative backward error for the linear systems solved by Ipopt during optimal power flow analysis for a 73-bus grid. The optimization problem has 1071 variables, 655 equality constraints, and 240 inequality constraints. Each  matrix is 2,237x2,237 with 11,297 nonzeros. The dashed line denotes the reciprocal of floating-point double precision.}
    \label{fig:solver_iterations}
\end{figure}

A typical sequence of linear systems generated by Ipopt is depicted in Figure \ref{fig:solver_iterations}, where the matrix condition number is superposed with the component-wise relative backward error computed with STRUMPACK (see Table~\ref{tab:errors} for the exact formula) for each of 84 matrix updates during the optimization solver run. Note that four of the first six iterations produce a matrix with condition number larger than $10^{20}$ (numerically singular). These correspond to failed iteration steps of the linesearch Newton method employed by Ipopt. The optimization solver can recover from this, for example, by adjusting the steplength. The condition number increases significantly somewhere after the 60th matrix update--- a normal feature of the continuation algorithm employed by Ipopt. It is at this stage that a robust linear solver is essential. All of our test cases are selected from the last few successful Ipopt iterations.

\subsection*{Literature review}

There is a vast body of research on benchmarking optimization solvers for their effectiveness and performance \cite{Dolan2002, Dolan2004}. However, not much attention has been devoted to the linear solvers employed within interior methods for optimization. The literature can be divided into three major categories: (a) papers that compare different optimization solvers without isolating the linear solver performance, including parameters and algorithm selection, (b) literature focusing on linear solver selection and performance, and (c) papers focusing on linear solvers for symmetric indefinite systems.

In papers from category (a), the CUTE, CUTEr, or CUTEst \cite{Gould2015} frameworks are often employed. For example, Schenk et al.~\cite{Schenk2007} applies a combinatorial (symmetric weighted matching based) pre-processing step to lower the condition number of the matrices arising in Ipopt, followed by a direct solver that uses supernodal static pivoting. The authors test this approach with CUTEr.

For category (b), authors are primarily concerned with finding suitable preconditioners \cite{Casacio2017, ONeal2005, Rees2007} for solving the linear systems iteratively. We cannot take advantage of these results, however, because the problems described herein are extremely ill-conditioned and thus require a direct linear solver, possibly with iterative refinement (IR), to achieve acceptable backward error levels.

For category (c), there is a long lineage of research available because the scientific community has been interested in efficient numerical linear solvers for (dense and sparse) symmetric indefinite systems dating back to the 1970s \cite{Bunch1971, Fletcher1976}. Bunch and Kaufman \cite{bunch1977some} introduced the pivoting strategy currently used in LAPACK for dense symmetric indefinite systems. This scheme was later improved in \cite{ashcraft1998accurate}, where the more stable \textit{bounded Bunch-Kaufman pivoting} is proposed. Duff et al.\ pioneered much of the field with a series of papers focused on symmetric indefinite systems \cite{duff1979direct,duff1983multifrontal,duff1991factorization,Duff2004,duff2005strategies}. Gould et al.\ \cite{gould2004numerical} reviewed the sparse symmetric solvers that are part of the mathematical software library HSL (formerly the Harwell Subroutine Library) and concluded that MA57 \cite{Duff2004} is currently the best general-purpose package in HSL. Hogg et al.\ \cite{hogg2013pivoting} improved the pivoting strategy of \cite{duff1991factorization} for challenging matrices. Duff et al.\ \cite{duff2007towards} focused on pivoting strategies and new approaches designed for parallel computers. Similarly, several recent papers \cite{hogg2014compressed,hogg2010indefinite,duff2018new} have considered improved algorithms targeting modern multicore processors and GPUs. Becker et al.~\cite{Becker2011} proposed a supernode-Bunch-Kaufman pivoting strategy that enables high-performance implementations on modern multicore processors. For a survey of the literature on sparse symmetric systems, see Davis et al.~\cite{davis2016survey}.  For a benchmark of several sparse symmetric direct solvers, see Gould et al.\ \cite{gould2007numerical}.


One of the few published resources, where different linear solver strategies are compared within Ipopt, is the PhD thesis of Jonathan Hogg \cite{Hogg2010}, who compares the performance of four linear solver frameworks: MA57, HSL\_MA77 \cite{Reid2009}, HSL\_MA79 \cite{Hogg2010b}, and Pardiso \cite{Pardiso}. Hogg devotes some consideration to the IR approach, and also analyzes parallel performance and scaling of the linear solver strategies, although when the thesis was published, the use of GPUs in scientific computing was still in its infancy. 

Our study is similar to the work of Tasseff et al.\ \cite{tasseff2019} and can be viewed as an extension. In that study, the authors analyze the performance of nine different linear solver frameworks (MA27, MA57, HSL\_MA77, HSL\_MA86, HSL\_MA97, MUMPS, PARDISO, SPRAL and WSMP) that can be applied in Ipopt to solve symmetric indefinite linear systems. All the linear solver packages taken into consideration by~\cite{tasseff2019} are based on a symmetric \LDLT{} factorization. One of them is SPRAL \cite{duff2018new}, which is also employed in our study. 

Several factors differentiate our study from \cite{tasseff2019}: (a) we consider linear solver software packages that use an unsymmetric LU factorization (which is potentially more stable), (b) we analyze GPU performance of the packages and compare the GPU performance to the CPU performance (the only GPU-accelerated software package tested in~\cite{tasseff2019} is SPRAL), (c) we focus on only one, not multiple test cases (they all come from the same application), (d) we extract matrices and right-hand sides from Ipopt instead of running the entire Ipopt code for each test case, and (e) in our comparisons, we not only compare performance but also check the solution accuracy. 

The major contribution of this paper is an in-depth comparison of well known linear solver packages for symmetric indefinite saddle-point systems arising from an interior method for optimization. 
The analysis includes CPU scaling (if applicable), GPU performance, and accuracy.

\section{Linear solver packages} \label{sec:linear_solvers}

We compare the performance and usability of five widely available linear solver packages: SuperLU~\cite{Li2005, SuperLUweb}, STRUMPACK~\cite{Rouet2016, SSTRUMPACKweb}, SPRAL/SSIDS~\cite{SPRAL}, cuSolver~ \cite{cuSolverweb}, PaStiX~\cite{PaStiX}, with MA57~\cite{Duff2004} as a reference. Of these packages, cuSolver is a proprietary product that can be obtained free of charge.  MA57 is licensed as part of the HSL library \cite{HSL}, while the others are open source software. We identified these solvers as best candidates to be used within the interior method context, based on results available in peer-reviewed literature.
We found that GPU support for PaStiX is less mature than for the other solvers and were not able to run all of our tests with it. PaStiX performance results are therefore not presented here, but partial results are included in the supplementary material at \url{https://github.com/kswirydo/linear-solver-review.
}

\subsection{SuperLU}

SuperLU is the most mature and established of the five packages evaluated. 
For the purpose of generating results presented in this study we used SuperLU\_dist, which is accelerated by MPI and CUDA. 

SuperLU is a direct linear solver for large, sparse unsymmetric systems. It employs Gaussian elimination (super-nodal) to compute the factorization
$$
P_rD_rAD_cP_c = LU,
$$
where $D_r$ and $D_c$ are diagonal matrices used to equilibrate (scale) the matrix $A$, and $P_c$, $P_r$ are column/row permutation matrices. 
If the solution is not accurate enough, it follows with an iterative refinement (IR) phase based on the Richardson iteration.

Unlike the serial version of SuperLU (which uses partial pivoting), the distributed version uses static pivoting based on graph matching \cite{Li1999}, which has better scaling properties \cite{Riedy2010}. The pivoting strategy determines the choice of $P_r$. Distributed SuperLU employs an $A^T+A$-based sparsity ordering, which determines $P_c$. SuperLU reorders both columns and rows; however, the user does not have the freedom to select a particular algorithm.

\subsection{STRUMPACK}

STRUMPACK (STRUctured Matrix PACKage was developed by Pieter Ghysels, Xiaoye Li, Yang Liu, Lisa Claus and other contributors at Lawrence Berkeley National Laboratory. STRUMPACK is intended to solve sparse and dense systems. STRUMPACK is targeted toward matrices that exhibit an underlying low-rank structure (such as block structure) that is exploited to construct the $L$ and $U$ factors. Unlike in SuperLU, the computed factors can be used as a preconditioner within an iterative linear solver. A multifrontal factorization is
employed.

STRUMPACK uses MC64 \cite{MA64} to permute and scale matrix columns (as does SuperLU). It follows the MC64 step with a nested dissection re-ordering~\cite{george1973}, which 
controls the amount of fill-in. The third pre-processing step is a reordering that helps reduce the ranks used for the HSS (hierarchically semi-separable) compression steps. After the pre-processing, there is a symbolic factorization phase to construct the elimination tree, and then actual multifrontal factorization. During the factorization, HSS is used to numerically compress the fronts. This strategy is typically only applied to large fronts (near the end of the factorization).  Depending on the settings and the matrix, the solve can be done through back substitution combined with either IR 
or an iterative solver (GMRES or BiCGStab). 

STRUMPACK is parallelized through MPI and CUDA. It uses cuBLAS and cuSolver, and relies on SLATE for LAPACK functions.

\subsection{cuSolver}

\verb|cuSolver| is a library provided by NVIDIA and distributed with CUDA Toolkit The functions provided by the library allow the user to solve a system (or multiple systems) of linear equations $Ax=b$, where $A$ does not have to be square. It computes a QR, LU, or \LDLT{} factorization of $A$, then performs upper and lower triangular solves. The library consists of three main components: \verb|cuSolverDN| (for dense matrices), \verb|cuSolverSP| (for sparse matrices) and \verb|cuSolverRF| (for matrix series in which the matrix values change but the nonzero pattern stays the same). The library comes with many helper functions and provides several algorithms for solving linear systems (using LU, QR, Cholesky, or least squares).

In this study, we tested \verb|cuSolverSP| because our matrices are sparse. However, the sparse interface appears to be less mature and to have limited functionality. For example, some of the functions are not implemented to run on GPUs, and \LDLT{} is available only through the dense interface. 

cuSolver is proprietary software and its source code is not available. The user supplies the matrix and right-hand side and the functions return a solution along with the error. The only parameter that can be chosen is the reordering with options: RCM (reverse Cuthill–McKee), AMD (approximate minimum degree), Metis, or none.

\subsection{SSIDS}

SSIDS is part of a large library known as SPRAL (Sparse Parallel Robust Algorithms Library). SPRAL was developed by the Computational Mathematics Group at Rutherford Appleton Laboratory in the UK.

Unlike SuperLU and STRUMPACK, SSIDS uses an \LDLT{} factorization to solve $Ax=b$ with $A$
symmetric but indefinite.  In a first pass, SSIDS computes the decomposition $A=PLD\left(PL\right)^T$, where $P$ is a permutation matrix, $L$ is unit lower triangular, and $D$ is either diagonal or block diagonal with $1 \times 1$ or $2\times 2$ blocks. A clear advantage of the \LDLT factorization is that it requires less storage, as only one triangular factor needs to be stored. SSIDS applies an \textit{a posteriori} threshold pivoting~\cite{duff2018new} to implement a scalable GPU pivoting strategy.

\subsection{Comparison}

The information concerning the basic functionality of the packages is collected in Table~\ref{tab:compbase}. 

\begin{table}[htb]  
\caption{\label{tab:compbase} Summary of linear solver package capabilities. The line labeled {\it Preconditioner} says if the packages compute an LU or \LDLT{} decomposition to be used as a preconditioner for an iterative linear solver. The capitalized names in the last row are commonly used names of BLAS routines.}

\medskip

\scalebox{0.6}{
\centering
\begin{tabular}{l|cp{2cm}p{1cm}p{1cm}c}
{\bf Solver }& SuperLU & STRUMPACK & SSIDS & cuSolver & PaStiX  \\\hline 
Exploits matrix symmetry & NO & NO & YES & NO & YES \\
Primary algorithm & LU & LU  & \LDLT{} &  LU, QR & \LDLT{}\\
Follows with IR & YES & YES & NO & NO & YES  \\
Preconditioner & NO & YES & NO & NO & YES \\
Supports reading RHS & NO & NO & NO & NO &NO \\
GPU-accelerated & YES & YES & YES & YES &YES\\
Multi-node & YES & YES & NO & NO & YES\\
Open source & YES & YES  & YES & NO &YES \\\hline
GPU-accelerated parts & GEMM &  GETRF, LASWP, GETRS, TRSM & SYRK, reord., fact., solve & QR fact., QR solve & fact. \\ \hline
\end{tabular}}
\end{table}

\section{Test case descriptions} \label{sec:tests}

All our test cases were generated by running 
optimal power flow analysis for synthetic power grid models developed by Texas A\&M University. The analysis was performed using Ipopt with MA57 as the linear solver.
We present results for three families of test cases named ACTIVSg  2000, 10k and 70k \cite{birchfield2016grid}. They closely resemble grids of Texas, US Western Interconnection and US Eastern Interconnection, respectively, and are representative of proprietary grid models used by power industry. Eastern Interconnection is the largest publicly available transmission grid model. Details on how these models are created are provided  in \cite{tamu}. Each family contains a set of linear problems (\ref{eq:mathproblemstatement}) from the last few Ipopt iterations.
The families of matrices and a more detailed description of the optimization problems are available in \cite{maack2020matrices}. The test case family characteristics (matrix size, number of matrices in the family, average number of nonzeros) are given in Table~\ref{tab:descr}.

\begin{table}[t]  
\caption{\label{tab:descr} Characteristics of the problems evaluated. Because the matrices are symmetric, only the triangular part of a matrix is stored in a Matrix Market file. {\it NNZ un} shows the number of nonzeros after unpacking (forming the full matrix with all of its nonzeros). The last column shows the number of matrices in each test case.}

\medskip

\centering
{\scriptsize
\begin{tabular}{ r |l | r r l l r  } 
\rowcolor{gray!40} &  & \multicolumn{5}{c}{Matrix characteristics} \\
 \rowcolor{gray!40} Model & Grid  & Rows & Cols & NNZ & NNZ un & \#   \\ 
\rowcolor{green!10} 2000 & ACTIVSg2000 & 56k   &56k  & 150k  &268k &3 \\
\rowcolor{green!10} 10k  & ACTIVSg10k & 238k &238k  &626k.8 &1,112k &5 \\
\rowcolor{green!10} 70k  & ACTIVSg70k & 1,640k  &1,640k   &4,320k &7,672k&5
\end{tabular}}
\end{table}
The results for other test cases are available in the repository: \url{https://github.com/kswirydo/linear-solver-review}.

\section{Test results} \label{sec:results}

All performance results presented in this section were obtained on one node of  Summit, the OLCF (Oak Ridge Leadership Computing Facility) supercomputer~\cite{OLCFSummit}.
Each node of Summit is equipped with two Power9 CPUs, each having 21 cores and six NVIDIA V100 GPUs. In all the test cases, we either used one rank and one GPU or multiple ranks with each rank having its own GPU. 

For each of the linear solver packages, we first present the results in terms of achieved (backward) error, followed by a short discussion. 
Each of the solver packages measures the solution accuracy differently. Table~\ref{tab:errors} collects the error descriptions. Afterwards, we provide the performance results.

\begin{figure}[tb]   
\centering
\includegraphics[width=0.45\textwidth]{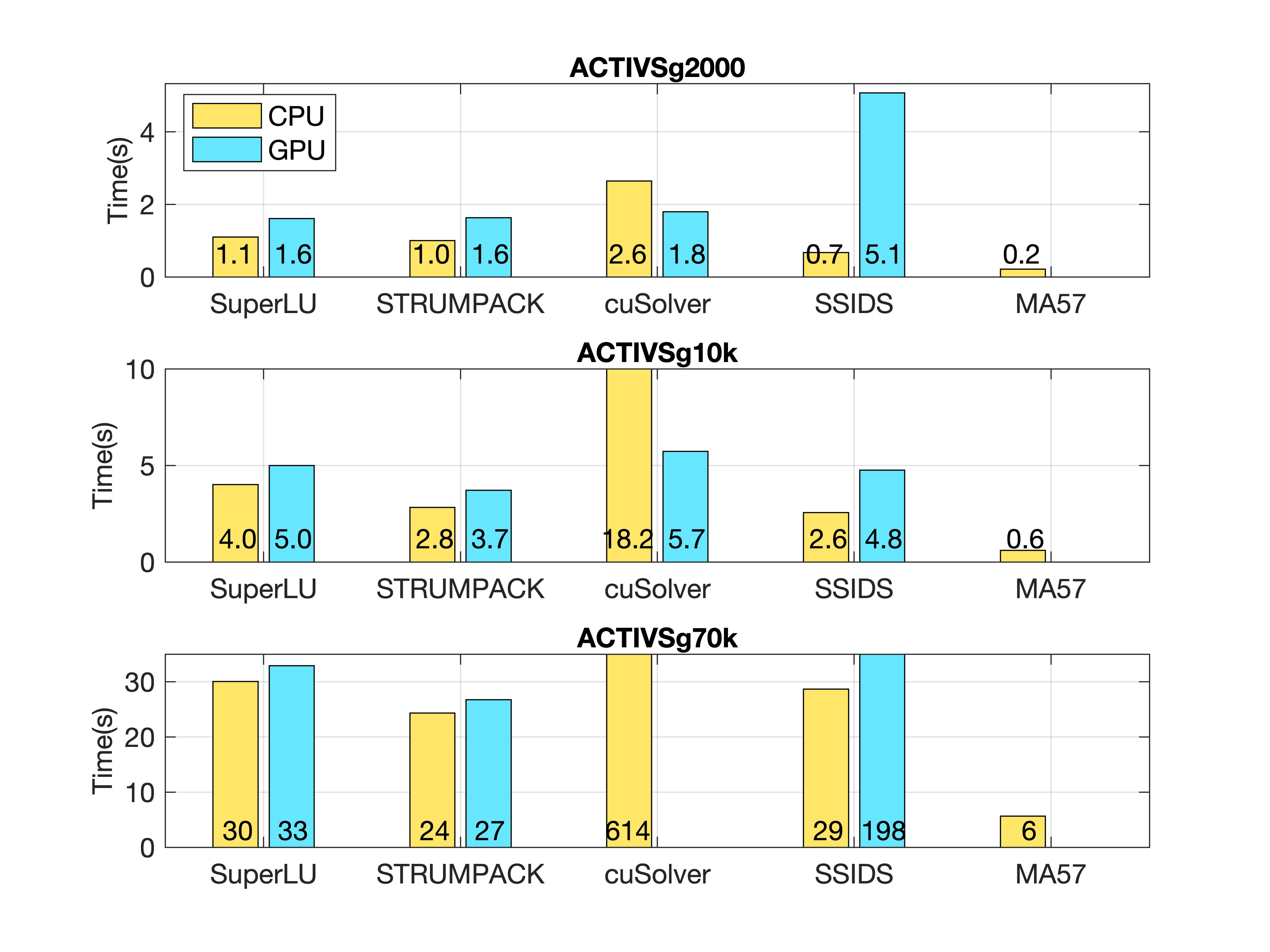}
\caption{\label{fig:comp} Comparison of the performance of all the solver packages tested. Times shown on the bars do not include matrix I/O nor error computation. Left: times without GPU acceleration. Right: times with GPU acceleration. 
}
\end{figure}

The results for all the tested packages are shown in Figure~\ref{fig:comp} together with results obtained using MA57, which is considered a standard solver for the type of linear systems discussed in this manuscript. In the figure, only best results obtained using  each of the packages are shown.

The results in the Figure are averaged per system. The systems for which the solver failed are  excluded from the averages.

\begin{table}[t]  
    \caption{\label{tab:errors} Different errors computed by the solver packages to measure solution accuracy. By $x_{\text{true}}$ we understand a known, true solution used in SuperLU testing.}
    
\medskip

\scriptsize
    \centering
    \begin{tabular}{m{2cm}m{3cm}m{2cm}}
        \rowcolor{gray!10} Name & Error formula & Used in \\\toprule
        Relative forward error (RFE) & $\displaystyle\frac{\|x - x_{\text{exact}}\|_{\infty}}{\|x\|_{\infty}}$ & SuperLU \\ \midrule
        Normwise Backward Residual Error (NBRE) v1 & $\displaystyle\max_i \frac{|r_i|}{s_i}\, \text{, where} \, s =|A||x|+|b|$ & SuperLU \\ \midrule
        Componentwise Relative Backward Error (CRBE) & $\displaystyle\max_{i=1,\ldots, n} \frac{\left|b_i - \sum_{j=1}^{n}a_{ij}x_j\right|}{|b_i| + \sum_{j=1}^{n}|a_{ij}| |x_j|}$ & STRUMPACK \\\midrule
        Relative residual & $\displaystyle \frac{\|b-Ax\|}{\|b\|}$ & STRUMPACK, cuSolver, SSIDS, PaStiX \\\midrule
        Normwise Backward Residual Error (NBRE) v2& $\displaystyle \frac{\|b-Ax\|_\infty}{\||b|+|A||x|\|_\infty}$ &  cuSolver, SSIDS,  PaStiX 
    \end{tabular}
\end{table}

\subsection{SuperLU}

\subsubsection{Convergence and the quality of solutions}

The SuperLU driver provided with the library computes NBRE v1 and RFE (see Table~\ref{tab:errors}). Obtaining the second quantity is not possible when the exact solution is not known {\em a priori}. Therefore, two sets of tests have been performed. In the first set, the right-hand side was generated by SuperLU based on an exact solution vector consisting of $1$s. This provides important information about the factorization quality (we call this {\it default RHS}). The second set was performed using right-hand sides matching the matrices (extracted from Ipopt) and read from files. In the second case, the reported error is the NBRE v1  achieved after IR. 

If a zero pivot is encountered, SuperLU stops the LU factorization. In such a case, if the exact solution is available, the error equals $1$ (vector of computed solutions is initialized to $b$) and no IR is performed.

The results are collected in Table~\ref{tab:superlu_res}, where:
\begin{itemize}
\item[] {\it NNZ inc.}\ is the ratio of the nonzeros in the $L$ and $U$ factors to the nonzeros in the original matrix;
\item[] {\it Av.\ RFE} is the average relative forward error (average through all matrices in each test case family);
\item[] {\it Av.\ NBRE} is the average norm-wise relative backward error
(reported after the last step of IR).
\end{itemize}

The matrices for which SuperLU failed were excluded from the averages.

\begin{table}[htb]  
\caption{\label{tab:superlu_res} The error levels achieved by SuperLU\_dist. {\it NNZ inc}: increase in nonzeros in $L$ and $U$ in relation to nonzeros in unpacked matrix; {\it Av.\ RFE}: average relative forward error (average through all matrices in each test case family); {\it Av.\ NBRE}: average NBRE v1.\ (error reported after last step of IR).} 

\medskip

\scriptsize
\centering
\begin{tabular}{ l |  r | l l } 
\rowcolor{gray!20} Name & NNZ inc. &av. RFE & av. NBRE \\
 ACTIVSg2000 &7.12x & $5.38 \cdot 10^{-7}$  &$2.73 \cdot 10^{-15}$ \\
 ACTIVSg10k  &6.52x & 2.74 $\cdot 10^{-8}$  & $5.78 \cdot 10^{-6}  $\\
  ACTIVSg70k &6.85x & $7.20$  $\cdot 10^{-3}$  & $1.045$
\end{tabular}
\end{table}

From Table~\ref{tab:superlu_res} one can infer that SuperLU produces solutions with sufficient accuracy for the optimization algorithm, even if the problem is close to being singular, and fails only on singular cases. We expect that these cases would be handled by the optimization algorithm controls, which adjust the iteration step and recompute the search direction.

\subsubsection{Performance results}

The computations performed by SuperLU can be divided into 8 phases: (a) equilibration (scaling), (b) row permutation, (c) column permutation, (d) symbolic factorization, (e) data distribution (through MPI), (f) (full) factorization, (g) solve, (h) IR.

\begin{figure*}[tbp]  
\scriptsize
\centering
\includegraphics[width=0.8\textwidth]{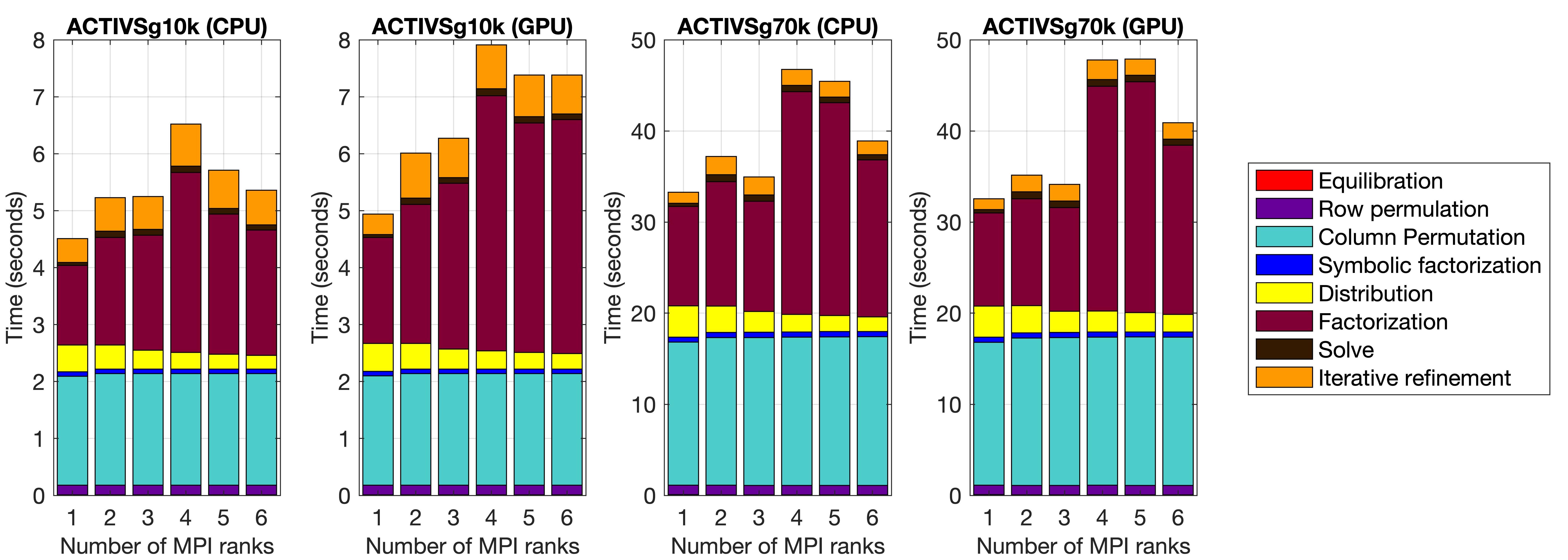}
\caption{\label{fig:superlu_defaultRHSCPU} SuperLU performance results using default right-hand sides. All times are given in seconds and all runs have one MPI rank per GPU. Left: results for ACTIVs10k, without and with GPU acceleration. Right: results for ACTIVSg70k without and with GPU acceleration.}
\end{figure*}

SuperLU can be compiled with or without GPU acceleration. Only phases (f) and (h) use acceleration. Figure~\ref{fig:superlu_defaultRHSCPU} and Table~\ref{tab:SuperLUperf} give results for both variants.
In this section, we only show results that were obtained using a fabricated RHS (generated by SuperLU). The extended set of results, including the results with correct RHS, are available in the repository.

\begin{table}[t]  
\caption{\label{tab:SuperLUperf}  SuperLU results for ACTIVSg2000 on one rank and one GPU. Eq = equilibration, RP = row permutation, CP = column permutation, SF = symbolic factorization, Dist = distribution, Fact = factorization, Ref = refinement. TOT is the total time (without I/O). The results show that only Fact and Ref use the GPU.}

\medskip

\centering
\scalebox{0.7}{
\begin{tabular}{ l | r r r r r r r r | r } 
\rowcolor{gray!40}\multicolumn{10}{c}{SuperLU}\\
\rowcolor{gray!20} & Eq & RP & CP & SF & Dist & Fact & Sol & Ref & TOT \\
\rowcolor{blue!20}CPU & 0.00 & 0.04 & 0.41  & 0.02 & 0.12 & 0.42  & 0.01 & 0.08 & 1.10 \\    
\rowcolor{blue!20} GPU &0.00 & 0.04 & 0.41  & 0.02 & 0.13 & 0.91  & 0.01 & 0.04 & 1.56 \\ 
\end{tabular}}
\end{table}

While analyzing Figure~\ref{fig:superlu_defaultRHSCPU}, we observe that certain parts of the linear solver scale quite well (e.g., the distribution for ACTIVSg70k decreases with the increase in the number of ranks). In the majority of cases, the highest computational cost is (as expected) the matrix factorization; however, for low MPI rank counts the computation is dominated by the column permutation step. A similar pattern can be observed for the GPU accelerated performance tests. The GPU times are much worse than  CPU times for the smaller cases and comparable to the CPU results for the two largest test cases (as expected from the matrix sizes). The reason for the time increase observed in the smallest test cases is due to data allocation/deallocation and movement. SuperLU uses a threshold value to decide whether the matrix is large enough to apply GPU acceleration. However, if the GPU accelerated version was used, SuperLU allocated and deallocated GPU memory regardless of the matrix size. 

\begin{figure}[tb]  
\centering
\includegraphics[width=0.5\textwidth]{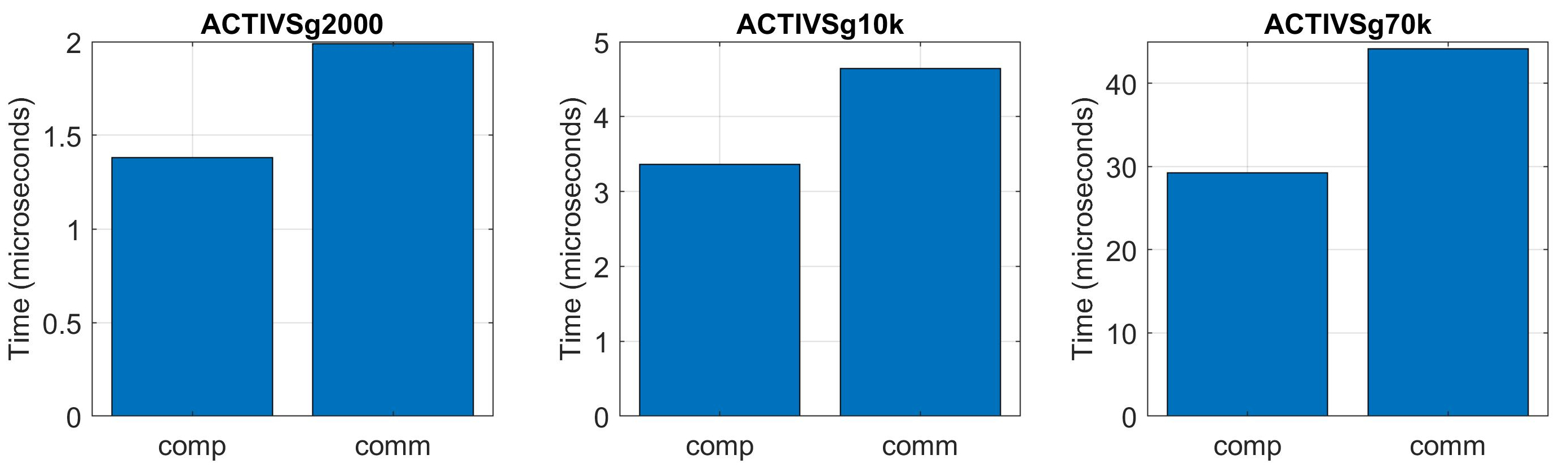}
\caption{\label{fig:superlu_GPU} Details of profiling SuperLU performance on one GPU. Only computation (comp) and communication (comm) times are shown. API calls took 722 ms (ACTIVSg2000), 861 ms (ACTIVSg10k), and 2230 ms (ACTIVSg70k).}
\end{figure}

To understand which part of the solve process occurs on the GPU, we launched compute jobs with \verb|nvprof| to obtain profiling results. 
We consider three categories of GPU operations: computations (computational kernels, either written by the SuperLU team or calling libraries such as cuBLAS), communication (explicit copies from Host to Device and from Device to Host) and API operations (cudaFree, cudaMalloc, cudaStreamCreate/Synchronize, cudaMallocHost, cudaMallocManaged, etc). All the detailed timing tests were performed on one GPU. 
The results are presented in Figure~\ref{fig:superlu_GPU}.

A detailed GPU profiling explains why the time increases significantly compared to the CPU-only case; the increase is largely due to GPU memory allocations, copies, and frees leading to a large startup cost.

We conclude that for all of our matrices, the GPU does not improve the overall performance. While it can be argued whether or not the API calls should be a part of the statistics (code analysis shows that these occur inside the linear solver functions), it is clear that the computations on the GPU are dominated by the cost of communication, as shown in Figure~\ref{fig:superlu_GPU}. Both are several orders of magnitude cheaper than the API operations (see the figure caption). 


\subsection{STRUMPACK}
\subsubsection{Convergence and quality of the solution}
\label{sec:strumpack_error}

STRUMPACK evaluates accuracy using component-wise relative backward error (CRBE) and relative residual; see Table~\ref{tab:errors}. The relative residual comes from using IR (Richardson by default).

We note that during the numerical phase of its multifrontal factorization, STRUMPACK pads matrices with zeros, which is not done in SuperLU and produces a factorization with more nonzeros compared to SuperLU (see~\cite{ghysels2016} for details of the factorization) .

For ACTIVSg2000, ACTIVSg10k, and ACTIVSg70k STRUMPACK with default settings (default depth of nested dissection (8), default iterative refinement options (Richardson)) fails because of zero pivots. For ACTIVSg2000, out of the three runs, one converged on the second matrix in the series, and the remaining two failed. 
This may be due to the randomized sampling applied to create the HSS structure. Linear systems as ill-conditioned as our test cases are very sensitive to a poor sampling and this may be the reason  why some of the test cases from the same family fail while some produce correct solutions.

 The lack of convergence can be prevented by adjusting the STRUMPACK parameters. It was determined that switching the IR algorithm to GMRES improved the relative residual.  However, we needed to increase the value of a parameter that determined the depth of nested dissection reordering (and resulting level of fill) to prevent zero pivots.  For ACTIVSg2000, this parameter was set to $25$; for ACTIVSg10k it was set to $50$; and for ACTIVSg70k it was set to $64$.

The numerical errors, summarized in Table~\ref{tab:strumpack_error}, are generally comparable with SuperLU.

\begin{table}[htb]  
\caption{\label{tab:strumpack_error} The error levels achieved in STRUMPACK. The results for ACTIVSg test cases were achieved by increasing the nested dissection parameter to $25$, $50$, $64$, respectively. All numbers computed using the default RHS.}

\medskip

\scriptsize
\begin{center}
\begin{tabular}{ l |r | l l } 
\rowcolor{gray!20} Name & NNZ inc & av. RR & av.  CRBE  \\ 
  ACTIVSg2000 &22.30x & $5.26\cdot 10^{-10}$ &$2.05 \cdot 10^{-11}$ \\
 ACTIVSg10k &36.42xx & $2.48\cdot 10^{-3}$  & $1.07 \cdot 10^{-12}$ \\
ACTIVSg70k &46.29x & $7.19 \cdot 10^{-3}$  & $7.50 \cdot 10^{-11}$
\end{tabular}
\end{center}
\end{table}

\subsubsection{Performance results}

STRUMPACK computations can be split into $5$ phases: (a)  re-ordering, (b) symmetrization, (c) symbolic factorization, (d) factorization, (e) solve.
STRUMPACK can be compiled with and without GPU acceleration. 


Because the end-to-end compute time increased significantly with the number of MPI ranks, we only present results for one rank (using four cores). A full set of results, including scaling up to $42$ ranks, can be found in the repository. When we examine the timings for particular components, the scaling properties are much better. For example, the solve time improves even for very small test cases, and the factorization time decreases for larger test cases. Results using a single rank are in Table~\ref{tab:STRUMPACK}.




\begin{table}[t]   
\caption{\label{tab:STRUMPACK}  STRUMPACK  results for all the test cases. All computations were performed on a single CPU (and single GPU). The nested dissection parameter was increased for the 3 test cases to 25, 50, and 64. Richardson IR was used for ACTIVSg10k on the GPU. For the remaining results, GMRES IR was used. ND=nested dissection, Sym=symmetrization, SF=symbolic factorization, Fact=numeric factorization, Sol=solve, TOT=total time (without I/O).}

\medskip

\centering
\scriptsize
\begin{tabular}{ l | r r r r r | r } 
\rowcolor{gray!40}\multicolumn{7}{c}{STRUMPACK}\\


\rowcolor{gray!20} &ND & Sym & SF & Fact & Sol & TOT \\ 
\rowcolor{blue!20} ACTIVSg2000 (CPU) & 0.45  &  0.08  &  0.01  &  0.07 &  0.03  &  0.64\\
\rowcolor{blue!20} ACTIVSg2000 (GPU) & 0.45  &  0.07  &  0.01  &  0.84  &  0.02  &  1.39 \\
\rowcolor{gray!20} ACTIVSg10k (CPU) & 1.86  &  0.33  &  0.024  &  0.50  &  0.11  &  2.82 \\
\rowcolor{gray!20} ACTIVSg10k (GPU) &  1.86  &  0.33  &  0.025  &  1.41  &  0.08  &  3.71 \\
\rowcolor{blue!20} ACTIVSg70k (CPU) &  15.14  &  2.30  &  0.14  &  5.69  &  1.08  &  24.35 \\
\rowcolor{blue!20} ACTIVSg70k (GPU) &  15.15  &  2.34 &  0.15  &  7.93  &  1.19  &  26.76
\end{tabular}
\end{table}

\begin{figure}[t]   
\centering
\includegraphics[width=0.5\textwidth]{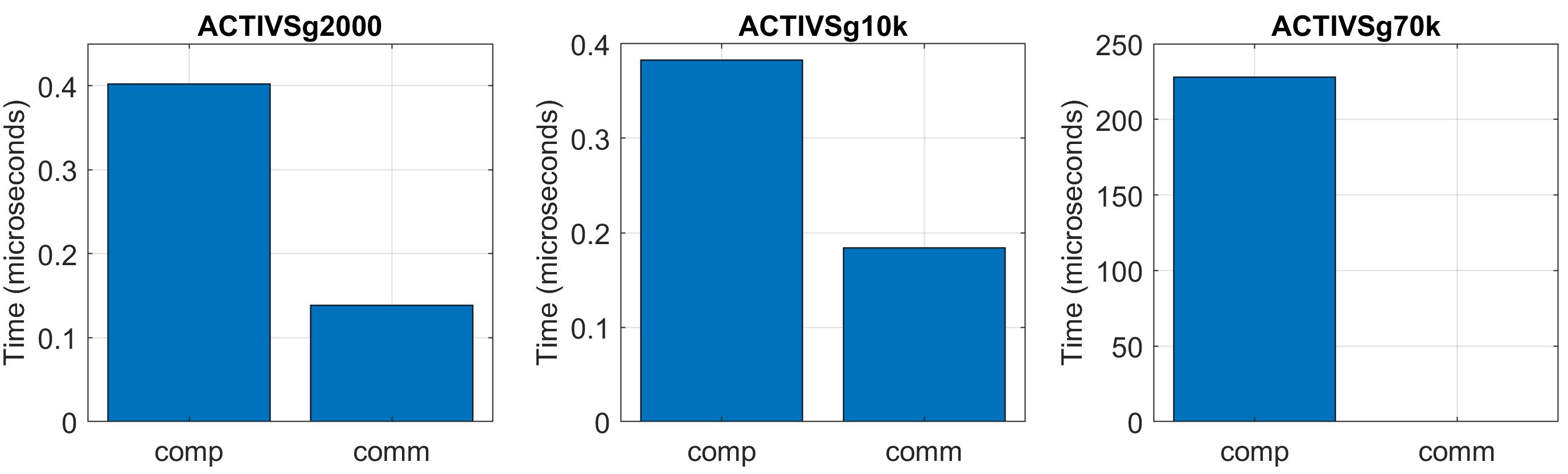}
\caption{\label{fig:STRUMPACK_GPU} Details of profiling STRUMPACK performance on the GPU. API calls are not shown; they took 1674 ms (ACTIVSg2000), 1649 ms (ACTIVSg10k), and 2234 ms (ACTIVSg70k). The discrepancy in the timing between the full stack timing as in Table \ref{tab:STRUMPACK} and the detailed GPU timing shown here originates in the difference between timing strategies. In the first case, \texttt{MPI\_Wtime()}  timer is used. In the second case, the code was run with \texttt{nvprof}, which adds overhead due to the nature of cuda profiling.}
\end{figure}

Next, we profile GPU performance. STRUMPACK uses a GPU if the matrix is large enough or more than one rank is specified. In order to observe any computations on the GPU, we need to use at least two MPI ranks or a large enough test case.  All results in Figure~\ref{fig:STRUMPACK_GPU} were obtained using two MPI ranks, except the last case which was tested with one MPI rank and one GPU.
See Figure~\ref{fig:STRUMPACK_GPU}.

For all problems except the last, the only computation that occurs on the GPU is the triangular solve. For these small cases, the factorization is not GPU-accelerated, but the memory allocation and memory copying that happen during the factorization are expensive and can explain the poor performance. For the last test case, many other computations (e.g., GEMM) are off-loaded to the GPU.

We conclude that there is no benefit to using GPU-accelerated STRUMPACK. In other cases, the overhead arising from the memory allocations outweighs any benefits. If a sequence of matrix equations is solved, and the GPU memory allocation (and free) happens only once, this cost is amortized over the optimization solver iterations. 

\subsection{cuSolver}

We tested linear solvers based on LU factorization and QR factorization. The linear solver function based on least squares was also appropriate but its performance was poor and therefore not included here. 

For a given matrix and right-hand side, the LU and QR linear solvers both return the solution $x$ and an integer that equals $-1$ if the matrix is nonsingular, and the index of the first zero row otherwise. 

All results in this section were computed using cuSolver from CUDA Toolkit version 10.1.243.

\subsubsection{Convergence and quality of solution}

The error is measured as relative residual, and as NBRE v2; see Table~\ref{tab:errors}. Note that for the last test case, the QR factorization on GPU fails. 
The results are shown in Table~\ref{tab:cuSolver_error}.

\begin{table}[t]   
\caption{\label{tab:cuSolver_error} The error levels achieved in cuSolver. Metis reordering was used for all the test cases. For ACTIVSg70k test case, QR method fails regardless of reordering.}

\medskip

\scriptsize
\begin{center}
\begin{tabular}{ l l | l l  } 
\rowcolor{gray!20}\multicolumn{4}{l}{ACTIVSg2000}  \\
&LU & $9.07\cdot 10^{-17}$& $6.10\cdot 10^{-23}$ \\
&QR & $2.96\cdot 10^{-16}$& $3.30\cdot 10^{-22}$ \\
\rowcolor{gray!20}\multicolumn{4}{l}{ ACTIVSg10k}  \\ 
&LU &  $1.41\cdot 10^{-10}$& $1.30\cdot 10^{-21}$ \\
&QR &  $3.47\cdot 10^{-11}$& $3.24\cdot 10^{-22}$ \\
\rowcolor{gray!20}\multicolumn{4}{l}{ ACTIVSg70k}  \\ 
&LU &  $9.98\cdot 10^{-7}$& $4.53\cdot 10^{-21}$ \\
&QR &  $--$& $--$ \\
\end{tabular}
\end{center}
\end{table}

\subsubsection{Performance results}

In addition to providing dimensions and CSR structure pointers, the user can choose the reordering (0: no reordering, 1: symrcm, 2: symamd, or 3: csrmetisnd) and the \textit{tolerance} to decide if singular or not.

While the tolerance has a rather marginal influence on both performance and solution quality, the type of reordering has a huge impact on the time. Thus, we set the tolerance to $10^{-12}$ for all tests and vary the reordering. 

As discussed earlier, both LU- and QR-based solves are {\it black box}. For a given matrix problem, they output the solution vector without information on the factors; thus we are unable to present detailed timing results as we did for SuperLU and STRUMPACK.

\begin{figure}[htb]   
\centering
\includegraphics[width=0.30\textwidth]{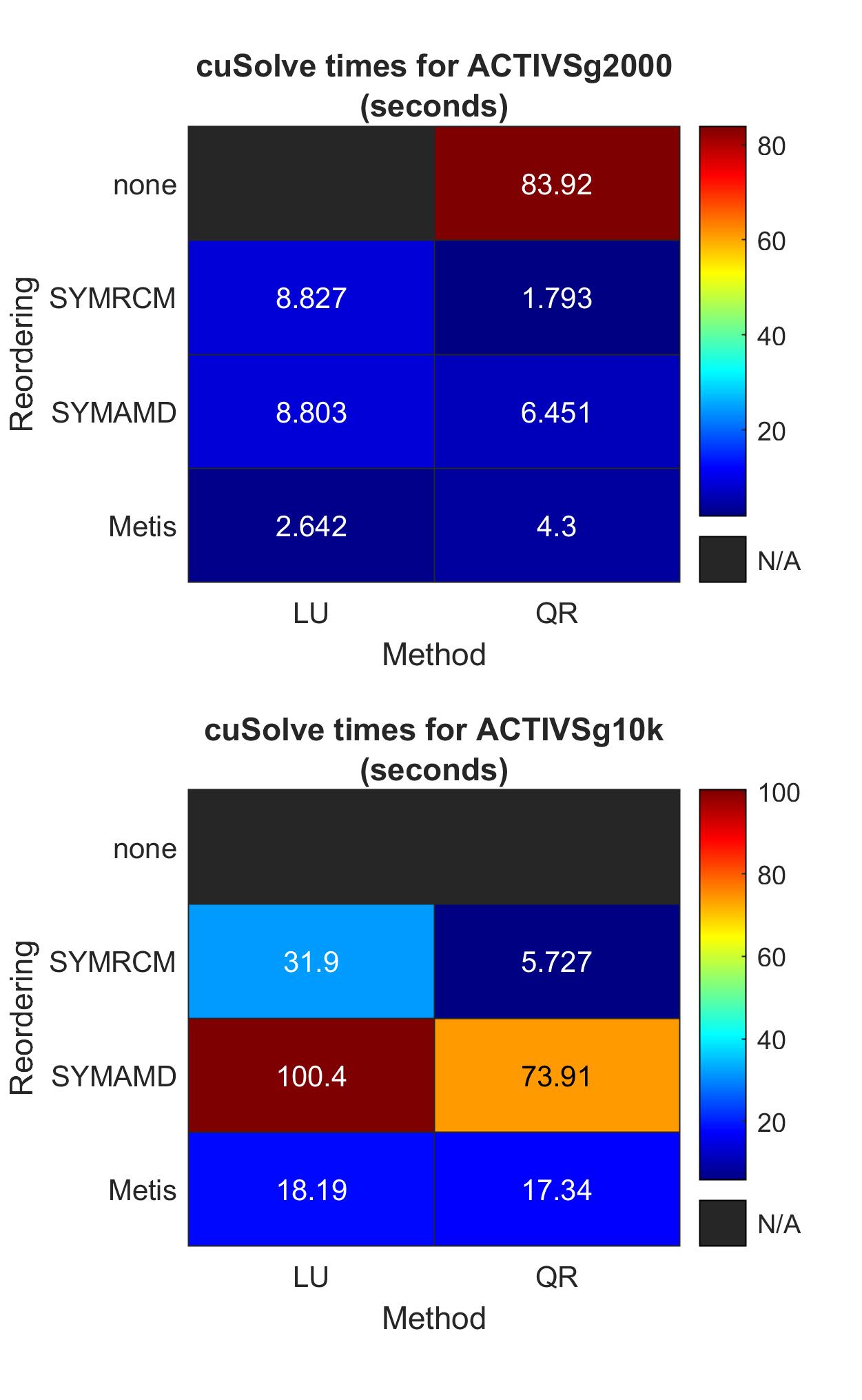}
\caption{\label{fig:cuSolver_perf}Performance results for cuSolver with different re-ordering options. Note that the ratio between timings for best and worst reordering is sometimes 5, and larger than 1000 for very small test cases.}
\end{figure}

While cuSolver is a GPU-accelerated library, only the QR solve function has a GPU interface. LU factorization is performed on the CPU. Hence the compute time for the LU solver was measured with a system timer and the compute time for the QR solver was measured with CUDA events. 
Results are shown in Figure~\ref{fig:cuSolver_perf}.

For ACTIVSg2000, CPU-based LU is faster than GPU-based QR, and again reordering has a huge impact on the time. Since even for the relatively small ACTIVSg2000, cuSolver was not able to complete the computation in 100 minutes on Summit, we did not try to run ACTIVSg10k without reordering. The only reordering that works for the ACTIVSg70k case is $3$ (metis). The time is $\approx 614$s per system (for LU) and the QR factorization fails regardless of the setting. The LU factorization is much slower than STRUMPACK. For the 70k case, the LU time in cuSolver can be lowered by a factor of $3.5\times$ if the matrix is scaled using Ruiz scaling~\cite{ruiz2001} prior to solving; however, STRUMPACK is still faster. The best results for the two largest cases are shown in Figure~\ref{fig:comp}.

The results in Figure~\ref{fig:cuSolver_perf} display the cost of solving the systems, but do not include the cost of data allocation, matrix read, and error evaluation. These costs were (averaged per system): 1.86s for ACTIVSg2000, 2.67s for ACTIVSg10k, and 9.97s for ACTIVSg70k.

\subsection{SSIDS}

\subsubsection{Convergence and quality of solution}

The error is measured in the same way as for cuSolver, as the scaled residual and NBRE v2; see Table~\ref{tab:errors}.
 The results are summarized in Table \ref{tab:SSIDS_error}:

\begin{table}[htb]   
\caption{\label{tab:SSIDS_error} The error levels achieved in SSIDS. }

\smallskip

\scriptsize
\begin{center}
\begin{tabular}{ l l | l l  } 
\rowcolor{gray!40} \multicolumn{2}{l|}{Name} & av. RR & av.  NRBE \\ 
\rowcolor{gray!20}\multicolumn{4}{l}{ACTIVSg2000}  \\
&CPU & $1.85\cdot 10^{+6}$& $2.35\cdot 10^{-12}$ \\
&GPU & $2.51\cdot 10^{-5}$& $2.25\cdot 10^{-12}$ \\
\rowcolor{gray!20}\multicolumn{4}{l}{ ACTIVSg10k}  \\ 
&CPU &  $1.67\cdot 10^{+4}$& $7.98\cdot 10^{-11}$\\
&GPU &  $5.99\cdot 10^{-1}$& $2.84\cdot 10^{-10}$ \\
\rowcolor{gray!20}\multicolumn{4}{l}{ ACTIVSg70k}  \\ 
&CPU &  $3.67\cdot 10^{+6}$& $2.55\cdot 10^{-12}$ \\
&GPU &  $1.86\cdot 10^{0}$& $2.69\cdot 10^{-12}$ \\
\end{tabular}
\end{center}
\end{table}



{\em ACTIVSg2000}.
For the 2000 series, in most cases GPU-accelerated SSIDS resulted in a segmentation fault. However, if run multiple times, it would eventually produce results, suggesting a possible write-past-the-end-of-array bug. For this series of problems, SSIDS exhibits good performance as long as no segmentation fault occurs. 

{\em ACTIVSg10k}.
While SSIDS is generally fast for these matrices, the relative residuals remain large. However, the norm-wise relative backward error is consistently very low.

{\em ACTIVSg70k}.
Only four of the five systems can be solved with the GPU-accelerated version. The CPU version again produces solutions with high scaled residuals and low relative backward error.

We conclude that SSIDS works very well for small problems (see the full set of results in the repository) It produces solutions with very small norm-wise relative backward errors.. For the RTS test cases, SSIDS is able to solve systems that none of the other 
packages can solve, which is impressive considering the lack of IR. After carefully studying the available documentation, we did not find a good explanation for the discrepancy between the relative residual values achieved on the CPU and on the GPU SSIDS.

\subsubsection{Performance results}

\begin{figure}[t]   
\centering
    \includegraphics[width=0.5\textwidth]{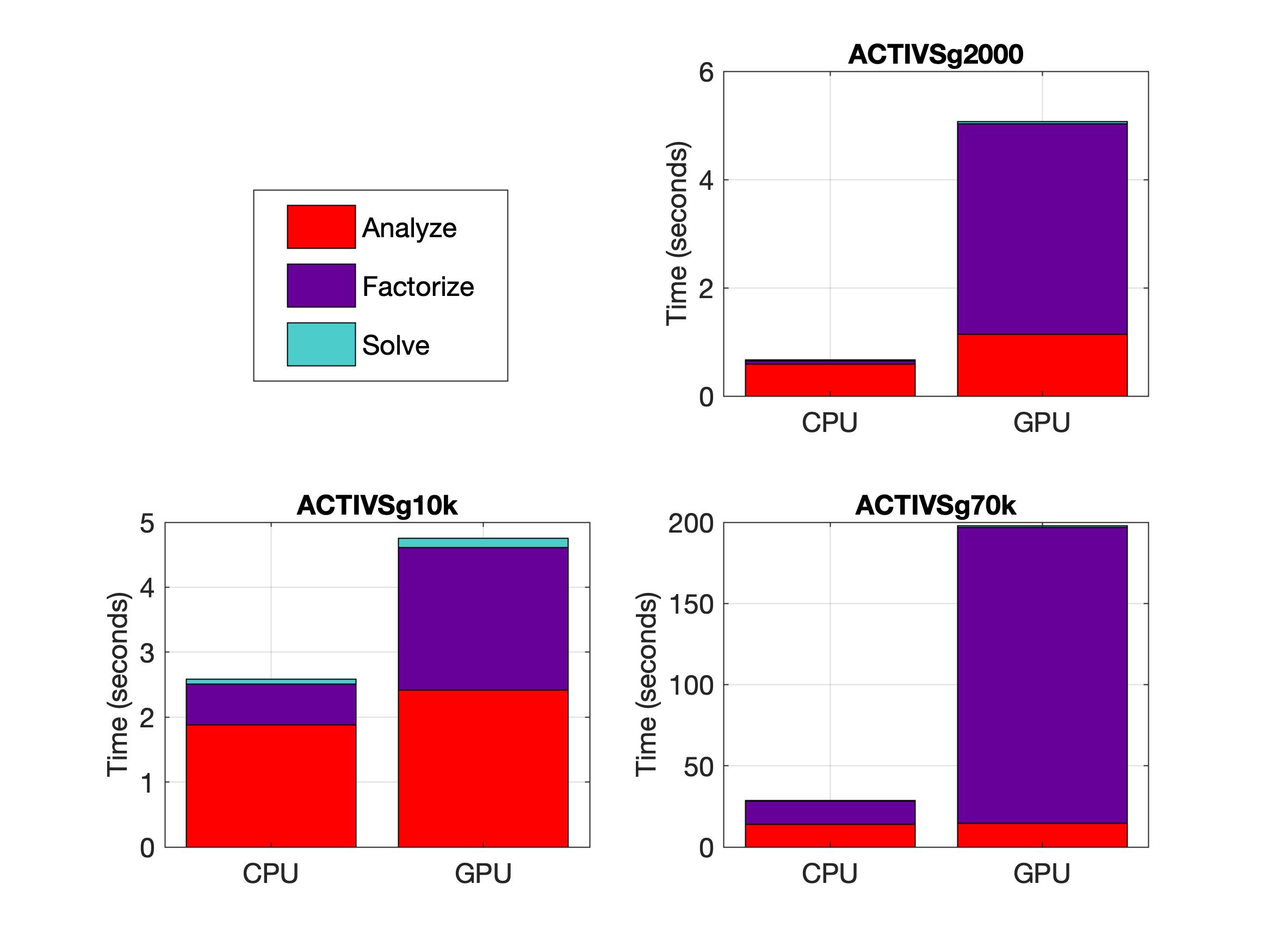}
    \caption{SSIDS performance on CPU and GPU. \label{fig:SSIDS}}
\end{figure}

\begin{figure}[tb]  
\centering
\includegraphics[width=0.5\textwidth]{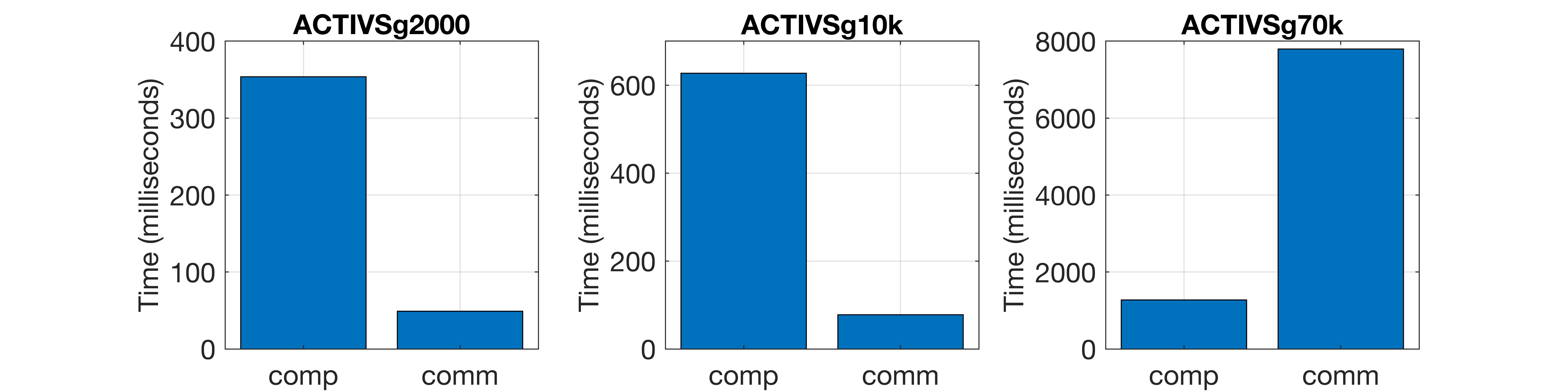}
\caption{\label{fig:ssids_GPU_prof} Details of profiling SSIDS performance on one GPU. Only computation (comp) and communication (comm) times are shown. API calls took $\approx 11$ s (ACTIVSg2000), $\approx 10$s (ACTIVSg10k), and $\approx 35$s (ACTIVSg70k).}
\end{figure}

SSIDS computations can be split into three basic phases: (a) analysis, (b) factorization, and (c) solve; this is how we time SSIDS. The results in Figure~\ref{fig:SSIDS} indicate that the GPU-accelerated version is much slower than the CPU-only version.
In comparison to STRUMPACK and SuperLU, SSIDS offloads more work to the GPU. However, because of the nature of the systems solved, the created {\it fronts} are very small. To factor them on the GPU, SSIDS allocates separate piece of memory at every level, which explains the huge increase in the time cost of the factorization compared to the CPU version. Profiling details are given in Figure~\ref{fig:ssids_GPU_prof}.

Random segmentation faults suggest the package is not quite mature.

\section{Comparison} \label{sec:summary}

The performance comparison of tested linear solvers is shown in Figure~\ref{fig:comp}. The results displayed do not include I/O and error evaluation.

CPU results: For the smallest test cases (available in the repository), cuSolver leaves other packages behind, but as the matrix size grows, the situation changes. SSIDS is the fastest for ACTIVSg2000; as the matrix size grows, STRUMPACK becomes the fastest.

GPU results:  As for the CPU, cuSolver can solve the smallest test case linear systems the fastest, but other packages become faster as the matrix sizes increases. SuperLU, STRUMPACK and cuSolver take roughly the same time for ACTIVSg2000, and STRUMPACK becomes a winner for the two largest test cases.

The most robust linear solver, able to process singular or close-to-singular systems, is SSIDS. 

All the  linear solver packages exhibit better CPU than GPU performance, which can be possibly mitigated if the memory used to store the matrices and factors is allocated once and reused over and over.  None of the 
packages scales well in distributed environments, which is expected from
our problem sizes. However, individual components of the solution process---such as solve and factorization---scale very well for some of the linear solvers (e.g., STRUMPACK).

\section{Conclusions and future work} \label{sec:conclusion}

In this study we tested five linear solver packages on challenging test problems arising from optimal power flow analysis for power grids. The test cases are linear problems (\ref{eq:mathproblemstatement}) that an interior-point optimization method hands off to the linear solver.

A common observation for the linear solver software is the lack of parallel scalability. If we look at the cost of the solution process (without I/O), the time decreases with the number of ranks, but only for certain components of the solution process. This lack of parallel scaling  is due to the properties of the matrices, which stress-test the linear solver packages beyond what they were designed for. The authors of all codes used in this paper were contacted to ensure the tests were run correctly. 

Based on the data collected for factorization times, we conjecture that the 
pivot strategies in LU and \LDLT{} solvers are among the main factors that explain the GPU performance. This is due to the 
irregular memory accesses and data movement associated with pivoting. The computations are memory bound: there is not enough computation per data movement, as can be inferred from our GPU profiling results for SuperLU and STRUMPACK.

In the case of SuperLU, STRUMPACK and SSIDS, there are excessive GPU memory allocations associated with solving multiple systems, leading to larger than necessary factorization times being reported. Therefore it is difficult to separate the cost of memory allocations from the effects of pivoting on performance for these 
solvers. For SSIDS, we observe that factorization times increase by $10\times$ over the CPU times; profiling results strongly suggests this happens because of repeated device memory allocations.

The two largely robust linear solver packages are STRUMPACK and SSIDS. To solve highly ill-conditioned systems, STRUMPACK might need to have its parameters adjusted to allow for more fill, but it is capable of producing high-quality solutions that have both low norm-wise relative backward error and low residual error. SSIDS is able to solve some of the systems that cause all other packages to fail. The advantage of SSIDS is that it respects matrix symmetry, and hence has lower storage requirements.

The cuSolver from NVIDIA consistently produces the lowest norm-wise relative backward error (NBRE v2) across all the test problems and thus represents the most accurate suite of algorithms evaluated in our study. The cuSolver for GPU is based on QR factorization, which in general exhibits better backward stability properties than LU 
\cite{Paige90}.

Future work will focus on FGMRES IR \cite{MarioA, MarioB}, with the possible application of mixed precision \cite{Erin} for increased execution speed on the GPU. To accelerate the solver phase we plan to explore iterative solution of sparse triangular systems as proposed by Chow et al.\ \cite{Edmond}.

\section*{Acknowledgements}
This research was supported by the Exascale Computing Project (17-SC-20-SC), a collaborative effort of the U.S. Department of Energy Office of Science and the National Nuclear Security Administration. We thank Chris Oehmen and Lori Ross O'Neil for critical reading of the manuscript and helpful suggestions. We also thank Tim Carlson and Kurt Glaesmann of Research Computing for their support, as well as Kestor Gokcen and Jiajia Li for providing utilities for matrix format conversion (all from Pacific Northwest National Laboratory). Finally, we acknowledge the support from the Oak Ridge Leadership Computing Facility, in particular a great deal of help from Philip Roth. 

\small
\bibliographystyle{plain}
\bibliography{references}
\end{document}